\documentclass[12pt]{article}
\usepackage{graphicx}
\usepackage{amsmath,amsthm,amssymb,enumerate}
\usepackage{euscript,mathrsfs}
\usepackage[left=2cm,right=2cm,top=3.5cm,bottom=3.5cm]{geometry}
\usepackage{color}
\catcode`\@=11 \@addtoreset{equation}{section}

\catcode`\@=12

\allowdisplaybreaks

\newtheorem{Theorem}{Theorem}[section]
\newtheorem{Proposition}[Theorem]{Proposition}
\newtheorem{Lemma}[Theorem]{Lemma}
\newtheorem{Corollary}[Theorem]{Corollary}

\theoremstyle{definition}
\newtheorem{Definition}[Theorem]{Definition}

\newtheorem{Remark}[Theorem]{Remark}

\newcommand{\bTheorem}[1]{
\begin{Theorem} \label{T#1} }
\newcommand{\eT}{\end{Theorem}}

\newcommand{\bProposition}[1]{
\begin{Proposition} \label{P#1}}
\newcommand{\eP}{\end{Proposition}}

\newcommand{\bLemma}[1]{
\begin{Lemma} \label{L#1} }
\newcommand{\eL}{\end{Lemma}}

\newcommand{\bCorollary}[1]{
\begin{Corollary} \label{C#1} }
\newcommand{\eC}{\end{Corollary}}

\newcommand{\bRemark}[1]{
\begin{Remark} \label{R#1} }
\newcommand{\eR}{\end{Remark}}

\newcommand{\bDefinition}[1]{
\begin{Definition} \label{D#1} }
\newcommand{\eD}{\end{Definition}}

\newcommand{\vcg}[1]{{\pmb #1}}

\newcommand{\bFormula}[1]{
\begin{equation} \label{#1}}
\newcommand{\eF}{\end{equation}}

\newcommand{\Ov}[1]{\overline{#1}}

\newcommand{\DC}{C^\infty_c}

\newcommand{\aleq}{\stackrel{<}{\sim}}
\newcommand{\ageq}{\stackrel{>}{\sim}}

\newcommand{\vr}{\varrho}
\newcommand{\vre}{\vr_\ep}

\newcommand{\vte}{\vt_\ep}
\newcommand{\vue}{\vu_\ep}

\newcommand{\vt}{\vartheta}
\newcommand{\vu}{\vc{u}}
\newcommand{\vm}{\vc{m}}
\newcommand{\Ee}{E_{\ep}}

\newcommand{\vq}{\vc{q}}

\newcommand{\vme}{\vm_\ep}
\newcommand{\vc}[1]{{\bf #1}}

\newcommand{\Div}{{\rm div}_x}
\newcommand{\Grad}{\nabla_x}

\newcommand{\dx}{\,{\rm d} {x}}

\newcommand{\dt}{\,{\rm d} t }

\newcommand{\intO}[1]{\int_{\Omega} #1 \ \dx}

\newcommand{\vv}{\vc{v}}

\newcommand{\ep}{\varepsilon}

\newcommand{\I}{\mathbb{I}}
\renewcommand{\S}{\mathbb{S}}

\definecolor{Cgrey}{rgb}{0.85,0.85,0.85}
\definecolor{Cblue}{rgb}{0.50,0.85,0.85}
\definecolor{Cred}{rgb}{1,0,0}
\definecolor{fancy}{rgb}{0.10,0.85,0.10}

\newcommand\Cbox[2]{%
    \newbox\contentbox%
    \newbox\bkgdbox%
    \setbox\contentbox\hbox to \hsize{%
        \vtop{
            \kern\columnsep
            \hbox to \hsize{%
                \kern\columnsep%
                \advance\hsize by -2\columnsep%
                \setlength{\textwidth}{\hsize}%
                \vbox{
                    \parskip=\baselineskip
                    \parindent=0bp
                    #2
                }%
                \kern\columnsep%
            }%
            \kern\columnsep%
        }%
    }%
    \setbox\bkgdbox\vbox{
        \color{#1}
        \hrule width  \wd\contentbox %
               height \ht\contentbox %
               depth  \dp\contentbox
        \color{black}
    }%
    \wd\bkgdbox=0bp%
    \vbox{\hbox to \hsize{\box\bkgdbox\box\contentbox}}%
    \vskip\baselineskip%
}

\date{}

\begin{document}


\title{Measure--valued solutions to the complete Euler system revisited}

\author{Jan B\v rezina \and Eduard Feireisl
\thanks{The research of E.F.~leading to these results has received funding from the
European Research Council under the European Union's Seventh
Framework Programme (FP7/2007-2013)/ ERC Grant Agreement
320078. The Institute of Mathematics of the Academy of Sciences of
the Czech Republic is supported by RVO:67985840.}
}

\date{\today}

\maketitle

\bigskip

\centerline{Tokyo Institute of Technology}

\centerline{ 2-12-1 Ookayama, Meguro-ku, Tokyo, 152-8550, Japan}

\bigskip

\centerline{Institute of Mathematics of the Academy of Sciences of the Czech Republic}

\centerline{\v Zitn\' a 25, CZ-115 67 Praha 1, Czech Republic}

\bigskip

\begin{abstract}

We consider the complete Euler system describing the time evolution of a general inviscid compressible fluid. We introduce a new concept of measure--valued solution based on the total energy balance and entropy inequality for the physical entropy without any renormalization. This class of so--called dissipative measure--valued solutions is large enough to include the vanishing dissipation limits of the Navier--Stokes--Fourier system. Our main result states that any sequence of weak solutions
to the Navier--Stokes--Fourier system with vanishing viscosity and heat conductivity coefficients
generates a dissipative measure-valued solution of the Euler system under some physically grounded constitutive relations.
Finally, we discuss the same asymptotic limit for the bi-velocity fluid model introduced by H.Brenner.

\end{abstract}



\section{Introduction}
\label{i}

We consider the \emph{complete Euler system} describing the time evolution of the mass density $\vr = \vr(t,x)$, the temperature $\vt = \vt(t,x)$
and the velocity field $\vu = \vu(t,x)$ of a compressible inviscid fluid:

\Cbox{Cgrey}{

\begin{eqnarray}
\label{i1} \partial_t \vr + \Div (\vr \vu) &=& 0, \\
\label{i2} \partial_t (\vr \vu) + \Div \left(\vr {\vu \otimes \vu}\right) + \Grad p  &=&0, \\
\label{i3} \partial_t \left( \frac{1}{2} \vr |\vu|^2 + \vr e \right) +
\Div \left[ \left( \frac{1}{2} \vr |\vu|^2 + \vr e + p \right) \vu \right] &=& 0.
\end{eqnarray}

}

\noindent
The system (\ref{i1}--\ref{i3}) contains the thermodynamic functions: The pressure $p(\vr, \vt)$ and the (specific) internal energy $e(\vr, \vt)$
depending on the state variables $\vr$, $\vt$ and satisfying Gibbs' relation
\begin{equation} \label{i4}
\vt D s = D e + p D \left( \frac{1}{\vr} \right).
\end{equation}
The new quantity $s$ appearing in (\ref{i4}) is the (specific) entropy. It follows from (\ref{i4}) that any smooth solution of (\ref{i1}--\ref{i3})
satisfies also the entropy balance
\begin{equation} \label{i5}
\partial_t (\vr s) + \Div (\vr s \vu) = 0
\ \mbox{or}\ \partial_t s + \vu \cdot \Grad s = 0.
\end{equation}
In the context of \emph{weak solutions}, the equation (\ref{i5}) is relaxed to the inequality
\begin{equation} \label{i6}
\partial_t (\vr s) + \Div (\vr s \vu) \geq 0,
\end{equation}
see e.g. Benzoni-Gavage, Serre \cite{BenSer}, Dafermos \cite{D4a}.
To avoid problems with physical boundaries, we restrict ourselves to the periodic boundary conditions, meaning the underlying physical domain $\Omega$
can be identified with the flat torus,
\[
\Omega = \left( [0,1] |_{\{ 0,1 \}} \right)^N, \ N=1,2,3.
\]
The problem is formally closed by prescribing the initial data
\begin{equation} \label{i15}
\vr(0, \cdot) = \vr_0, \ \vt(0, \cdot) = \vt_0,
\ \vu(0, \cdot) = \vu_0.
\end{equation}

In view of recent results based on the theory of convex integration, see \cite{FeKlKrMa}, weak solutions of (\ref{i1}--\ref{i3}), even if supplemented by
(\ref{i6}), are not uniquely determined by the initial data
as long as $N > 1$. As a matter of fact, for any \emph{piecewise constant} initial density $\vr_0$ and $\vt_0$, there exists
$\vu_0 \in L^\infty(\Omega; R^N)$, $N=2,3$ such that problem (\ref{i1}--\ref{i3}), (\ref{i6}), and (\ref{i15}) admits infinitely many weak (distributional) solutions on a given time interval $(0,T)$. This kind of result indicates that the measure--valued solutions that are supposed to capture possible oscillatory behavior of the weak solutions may be a suitable
concept for the Euler system.

In \cite{BreFei17}, we have introduced a concept of dissipative measure--valued (DMV) solution, based on postulating the total energy balance and a renormalized version of the entropy production equation, see Definition \ref{D2} below. In contrast with the standard approach used e.g. in a series of papers
by Fjordholm et al. \cite{FjKaMiTa}, \cite{FjMiTa2}, \cite{FjMiTa1}, based on the hypothetical $L^\infty-$bounds on the family of generating solutions, our definition covers
a more general class of objects that may be seen as suitable limits of weak solutions (cf. \cite{BreFei17}) or even certain numerical schemes.
In the present paper, we
introduce a slightly different definition of (DMV) solutions without entropy renormalization, see Definition \ref{D1} below. Our main result asserts that this kind of measure--valued solutions can be recovered as a vanishing viscosity limit of a natural physical approximation - the Navier--Stokes--Fourier system.
More specifically, we show
that any sequence of weak solutions to the Navier--Stokes--Fourier system, the existence of which is guaranteed by \cite{FeNo6},
generates a (DMV) solution of the
Euler system under some physically grounded constitutive relations. Finally, we discuss the
same asymptotic limit for the bi-velocity fluid model introduced by H.Brenner.

The paper is organized as follows. In Section \ref{PP} we introduce the basic concepts as well as the known results used in the text. In Section \ref{ZD},
we consider the vanishing dissipation limit of the weak solutions to the Navier--Stokes--Fourier system and show that they generate a (DMV) solution of the Euler system.
Finally, in Section \ref{BM}, we briefly discuss similar issues for a bi-velocity fluid model proposed by H.Brenner.

\section{Preliminary results}
\label{PP}

In this preliminary section, we introduce  conservative variables, the concept of (DMV) solution as well as other already known results used in the paper.

\subsection{Conservative variables}
\label{CV}

In certain situations, for instance in numerical simulations, it is more convenient to introduce the conservative variables:
The density $\vr$, the momentum $\vm = \vr \vu$ and the total energy
$E = \frac{1}{2} \vr |\vu|^2 + \vr e$ converting (\ref{i1}--\ref{i3}), (\ref{i6}) into

\Cbox{Cgrey}{

\begin{eqnarray}
\label{i7} \partial_t \vr + \Div \vm &=& 0, \\
\label{i8} \partial_t \vm + \Div \left(\frac{\vm \otimes \vm}{\vr} \right) + \Grad p  &=&0, \\
\label{i9} \partial_t E +
\Div \left[ \left( E +  p \right) \frac{\vm}{\vr} \right] &=& 0,\\
\label{i10}
\partial_t (\vr s) + \Div ( s \vm) &\geq& 0.
\end{eqnarray}

}

For the sake of simplicity, we focus on the case of polytropic gas, for which the pressure is related to the internal energy by the
caloric equation of state
\begin{equation} \label{i12bis}
p = (\gamma - 1) \vr e, \ \mbox{with the adiabatic exponent}\ \gamma > 1.
\end{equation}
Accordingly, we have
\begin{equation} \label{i12}
p = (\gamma -1) \left[ E - \frac{1}{2} \frac{|\vc{m}|^2}{\vr} \right]
\end{equation}
closing the system of equations (\ref{i7}--\ref{i9}). Under these circumstances, it is convenient to consider the specific entropy
$s = s(\vr, e)$ as a function of $\vr$, $e$, for which Gibbs' relation (\ref{i4}) yields
\[
\frac{\partial s}{\partial e}(\vr, e) = \frac{1}{\vt}, \ \frac{\partial s}{\partial \vr} (\vr, e) = - \frac{p}{\vt \vr^2},
\]
where the first relation can be seen as a definition of the absolute temperature. Moreover, it can be deduced from Gibbs' relation  (\ref{i4}) and  (\ref{i12bis})
that the entropy $s$ can be written in the form
\begin{equation} \label{entro}
s(\vr, e)  = S \left( \frac{p}{\vr^\gamma} \right) = S \left( \frac{(\gamma - 1) e }{\vr^{\gamma - 1}} \right)
\end{equation}
for a suitable function $S$.

\subsection{Thermodynamic stability}

In the original state variables $\vr$, $\vt$, the thermodynamic stability hypothesis reads
\begin{equation} \label{i16}
\frac{\partial p}{\partial \vr} (\vr, \vt) > 0,\
\frac{\partial e}{\partial \vt} (\vr, \vt) > 0 \ \mbox{for any}\ \vr, \vt > 0.
\end{equation}

In the conservative variables $(\vr, \vc{m}, E)$, this is equivalent to the statement that the total entropy
\begin{equation} \label{i16bisa}
\mathcal{S} (\vr, \vc{m}, E) = \vr S \left( \frac{(\gamma - 1)\left( E - \frac12\frac{|\vc{m}|^2}{\vr} \right) }{\vr^{\gamma}} \right)
\end{equation}
is a concave function of $(\vr, \vc{m}, E)$,
cf. Bechtel, Rooney, and Forest \cite{BEROFO}. It is a matter of direct computation to check
that, in terms of the function $S$ introduced in (\ref{entro}),
the condition (\ref{i16}) reduces to
\begin{equation} \label{i16bis}
(1 - \gamma) S'(Z) - \gamma S^{''}(Z) Z > 0 \ \mbox{for all} \ Z > 0.
\end{equation}
Note that the domain of definition of $S$ may not be $(0, \infty)$.
To see this, it is convenient to write $p$ and $e$ interrelated through
(\ref{i12bis}) as functions of $\vr$ and $\vt$. Accordingly, Gibbs' equation
(\ref{i4}) can be written in the form of Maxwell's relation
\[
\frac{\partial e}{\partial \vr } = \frac{1}{\vr^2} \left( p - \vt \frac{\partial p}{\partial \vt} \right) ,
\]
which, together with (\ref{i12bis}), gives rise to
\begin{equation} \label{i17}
p(\vr, \vt) = \frac{P(q)}{q^{\gamma}} \vr^\gamma \ \mbox{for a certain function} \ P, \ \mbox{where we have set} \ q = \frac{\vr}{\vt^{c_v}}, \ c_v = \frac{1}{\gamma - 1}.
\end{equation}
Furthermore, it follows from the second inequality in (\ref{i16}) that
\[
q \mapsto \frac{P(q)}{q^\gamma} \ \mbox{is a non--incresing function of}\ q,
\]
in particular,
\[
\lim_{\vt \to 0+} \frac{p(\vr, \vt)}{\vr^\gamma} = \Ov{p} \geq 0.
\]
We infer that the ``natural'' domain of definition of $S = S(Z)$ is the interval $(\Ov{p},\infty)$,
\[
S : (\Ov{p}, \infty) \to R.
\]
In addition, we define
\[
S(Z) = \left\{
\begin{array}{l}
- \infty \ \mbox{for}\ Z < \Ov{p},\\
\lim_{Z \to \Ov{p}+} S(Z) \in [-\infty, \infty) \ \mbox{for}\ Z = \Ov{p}.
\end{array}
\right.
\]
Accordingly, the total entropy $\mathcal{S} = \mathcal{S}(\vr, \vc{m}, E)$ is concave upper semi--continuous
in $[0, \infty) \times R^N \times [0, \infty)$ ranging in $[- \infty, \infty)$.

Similarly, we observe that the kinetic energy
\[
(\vr, \vc{m}) \mapsto \frac{1}{2} \frac{ |\vc{m}|^2 }{\vr} \ \mbox{is a convex function for}\ \vr > 0;
\]
whence we may define
\[
\frac{1}{2} \frac{ |\vc{m}|^2 }{\vr} = \left\{
\begin{array}{l}
0 \ \mbox{whenever} \ \vc{m} = 0, \ \vr \geq 0,\\
\infty \ \mbox{for}\ \vr = 0, \ \vc{m} \ne 0.
\end{array}
\right.
\]
The kinetic energy is therefore a convex lower semi--continuous function defined for $(\vr, \vm) \in [0, \infty) \times R^N$ ranging in
$[0, \infty]$.

\subsection{Relative energy}

The relative energy functional introduced in \cite{FeiNov10}, reads
\[
\begin{split}
\mathcal{E} \left(\vr, \vt, \vu \Big| \tilde{\vr}, \tilde{\vt}, \tilde{\vu} \right)
&= \frac{1}{2} \vr | \vu - \tilde{\vu} |^2 + H_{\tilde \vt}(\vr, \vt) -
\partial_\vr H_{\tilde \vt} (\tilde{\vr}, \tilde{\vt}) (\vr - \tilde \vr) - H_{\tilde \vt}(\tilde \vr, \tilde \vt)\\
H_{\tilde \vt}(\vr, \vt) &\equiv \vr \Big( e(\vr, \vt) - \tilde{\vt} s(\vr, \vt) \Big).
\end{split}
\]

Passing to the conservative variables,
\[
{\vm} = \vr {\vu}, \ {E} = \frac{1}{2} \vr | {\vu}|^2 + \vr e( \vr, \vt), \
\tilde{\vm} = \tilde \vr \tilde{\vu}, \ \tilde{E} = \frac{1}{2} \tilde \vr |\tilde{\vu}|^2 + \tilde \vr e(\tilde \vr, \tilde \vt),
\]
we can check by a bit tedious but straightforward manipulation that
\begin{equation} \label{i21}
\begin{split}
\mathcal{E} \left(\vr, \vt, \vu \Big| \tilde{\vr}, \tilde{\vt}, \tilde{\vu} \right) &\equiv
\mathcal{E} \left(\vr, E, \vm \Big| \tilde{\vr}, \tilde{E}, \tilde{\vm} \right)
\\
= - \tilde \vt &\Big[ \mathcal{S}(\vr, \vc{m}, E)  \\ &- \partial_{\vr} \mathcal{S}(\tilde \vr, \tilde {\vm}, \tilde E) (\vr - \tilde \vr)
- \nabla_{\vm} \mathcal{S}(\tilde \vr, \tilde {\vm}, \tilde E) \cdot (\vm - \tilde{\vm})
- \partial_{E} \mathcal{S}(\tilde \vr, \tilde {\vm}, \tilde E) (E - \tilde E)
  \\ &- \left.  \mathcal{S}(\tilde \vr, \tilde {\vc{m}}, \tilde E)             \right],
\end{split}
\end{equation}
where $\mathcal{S}$ is the total entropy introduced in (\ref{i16bisa}).
It is worth--noting that the expression in the brackets on the right--hand side of (\ref{i21}) coincides with the relative entropy \` a la Dafermos \cite{Daf4}.
In agreement with the discussion in the previous section, the relative energy plays a role of distance between $(\vr, \vt, \vu)$ and
$(\tilde \vr, \tilde \vt, \tilde{\vu})$, or, equivalently, between $(\vr, E, \vm)$ and $(\tilde{\vr}, \tilde{E}, \tilde{\vm})$,
as long as the thermodynamic stability hypothesis holds.

\subsection{Measure-valued solutions}
\label{MV}

In contrast with our preceding paper \cite{BreFei17}, we define a \emph{dissipative measure--valued (DMV) solution} with respect to the conservative variables
$(\vr, \vc{m}, E)$. Accordingly, the phase space for the associated Young measure is
\[
Q = \left\{ (\vr, \vc{m}, E) \ \Big| \ \vr \in [0, \infty), \ \vc{m} \in R^N, \ E \in [0, \infty) \right\}.
\]

\begin{Definition} \label{D1}
A \emph{dissipative measure--valued solution} to the problem (\ref{i7}--\ref{i10}) consists of a family of \emph{parameterized probability measures} $\{ Y_{t,x} \}_{t \in (0,T),
x \in \Omega}$ and a non--negative function $\mathcal{D} \in L^\infty(0,T)$ called \emph{dissipation defect} satisfying:
\begin{itemize}
\item $Y \in L^\infty_{{\rm weak-(*)}}((0,T) \times \Omega; \mathcal{P}(Q))$, where $\mathcal{P}(Q)$ denotes the set of probability measures on $Q$;
\item
\begin{equation} \label{MV1}
\int_0^\tau \intO{ \left[ \left< Y_{t,x}; \vr \right> \partial_t \varphi + \left< Y_{t,x}; \vm \right> \cdot \Grad \varphi \right] } \dt = \left[ \intO{ \left< Y_{t,x}; \vr \right> \varphi } \right]_{t = 0}^{t = \tau}
\end{equation}
for a.a. $\tau \in (0,T)$ and any $\varphi \in C^\infty([0,T] \times \Omega)$;
\item
\begin{equation} \label{MV2}
\begin{split}
\int_0^\tau &\intO{ \left[ \left< Y_{t,x}; \vm \right> \cdot \partial_t\vcg{\varphi} + \left< Y_{t,x} ; \frac{ \vm \otimes \vm}\vr \right> : \Grad \vcg{\varphi} + (\gamma - 1) \left< Y_{t,x}; \left(E - \frac12 \frac{|\vm|^2}\vr\right) \right> \Div \vcg{\varphi} \right] } \dt \\ &=
\left[  \intO{ \left< Y_{t,x}; \vc{m} \right> \cdot \vcg{\varphi} } \right]_{t = 0}^{t = \tau} + \int_0^\tau \int_{\Omega} \Grad \vcg{\varphi} : {\rm d}\mu_C
\end{split}
\end{equation}
for a.a. $\tau \in (0,T)$ and any $\vcg{\varphi} \in C^\infty([0,T] \times \Omega; R^N)$,
where $\mu_C$ is a (vectorial) signed measure on $[0,T] \times \Omega$;
\item
\begin{equation} \label{MV3}
\left[ \intO{ \left< Y_{t,x}; E \right> } \right]_{t = 0}^{t = \tau} + \mathcal{D}(\tau) = 0
\end{equation}
for a.a. $\tau \in (0,T)$;
\item
\begin{equation} \label{MV4}
\begin{split}
&\left[  \intO{  \left< Y_{t,x}; \vr s \left(\varrho, \frac{1}{\vr} \left( E - \frac{1}{2} \frac{|\vc{m}|^2}{\vr} \right) \right)  \right> \varphi } \right]_{t = 0}^{t = \tau}
\\ &\geq \int_0^\tau \intO{ \left[ \left< Y_{t,x} ; \vr s \left(\varrho, \frac{1}{\vr} \left( E - \frac{1}{2} \frac{|\vc{m}|^2}{\vr} \right) \right) \right> \right] \partial_t
\varphi} \dt \\ &+ \int_0^\tau \intO{ \left[
\left< Y_{t,x} ; s \left(\varrho, \frac{1}{\vr} \left( E - \frac{1}{2} \frac{|\vc{m}|^2}{\vr} \right) \right) \vc{m} \right> \cdot \Grad \varphi  \right] } \dt
\end{split}
\end{equation}
for a.a. $\tau \in (0,T)$, any $\varphi \in C^\infty([0,T] \times \Omega)$, $\varphi \geq 0$;
\item
\begin{equation} \label{MV5}
\left\| \mu_C \right\|_{\mathcal{M}([0, \tau] \times \Omega) } \leq c \int_0^\tau \mathcal{D}(t) \ \dt \ \mbox{for a.a.}\ \tau \in (0,T).
\end{equation}

\end{itemize}
\end{Definition}

\begin{Remark} \label{RM1}

The parameterized family of measures $\{ Y_{0,x} \}_{x \in \Omega}$ plays the role of initial conditions, cf. (\ref{i15}).

\end{Remark}

\begin{Remark} \label{RM2}

In Definition \ref{D1} we tacitly assume that all integrals are finite. In particular,
as
\[
\vr s \left(\varrho, \frac{1}{\vr} \left( E - \frac{1}{2} \frac{|\vc{m}|^2}{\vr} \right) \right) = \mathcal{S}(\vr, \vc{m}, E),
\]
where $\mathcal{S}$ is a concave function of all arguments, we have
\[
Y_{t,x} \left\{ (\vr, \vc{m}, E) \ \Big| \ \vr \geq 0, \ E - \frac{1}{2} \frac{|\vc{m}|^2}{\vr} \geq \frac{\Ov{p}}{\gamma - 1}\vr^\gamma  \right\} = 1
\]
for a.a. $(t,x) \in (0,T) \times \Omega$.

\end{Remark}

Following \cite{BreFei17}, we may replace (\ref{MV4}) by its renormalized version, namely
\begin{equation} \label{MV6}
\begin{split}
\Big[  &\intO{  \left< Y_{t,x}; \vr Z \left( s \left(\varrho, \frac{1}{\vr} \left( E - \frac{1}{2} \frac{|\vc{m}|^2}{\vr} \right) \right) \right) \right> \varphi } \Big]_{t = 0}^{t = \tau}
\\ &\geq \int_0^\tau \intO{ \left[ \left< Y_{t,x} ; \vr Z \left( s \left(\varrho, \frac{1}{\vr} \left( E - \frac{1}{2} \frac{|\vc{m}|^2}{\vr} \right) \right) \right) \right> \right] \partial_t
\varphi} \dt \\ &+ \int_0^\tau \intO{ \left[
\left< Y_{t,x} ; Z \left( s \left(\varrho, \frac{1}{\vr} \left( E - \frac{1}{2} \frac{|\vc{m}|^2}{\vr} \right) \right) \right) \vc{m} \right> \cdot \Grad \varphi  \right] } \dt,
\end{split}
\end{equation}
for a.a. $\tau \in (0,T)$, any $\varphi \in C^1([0,T] \times \Omega)$, $\varphi \geq 0$, where
$Z \in C(R)$ is a nondecreasing concave function and $Z(s) \leq Z_\infty$ for any $s \in R$.

\begin{Definition} \label{D2}

A \emph{renormalized (DMV) solution} to (\ref{i7}--\ref{i10}) consists of a family of parameterized probability measures $\{ Y_{t,x} \}_{t \in (0,T),
x \in \Omega}$ and a non--negative function $\mathcal{D} \in L^\infty(0,T)$ that satisfy all the requirements of Definition \ref{D1} except (\ref{MV4}), which is
replaced by (\ref{MV6}).

\end{Definition}

The concept of renormalized (DMV) solution is motivated by a similar definition of the weak solutions introduced in Chen and Frid \cite{CheFr2}. As shown in
\cite[Section 2.1.1]{BreFei17}, the renormalized solutions enjoy certain minimum principle, in particular
\[
Y_{0,x} (\left\{ s(\vr, \vc{m},E) \geq s_0 \right\}) = 1 \ \mbox{implies}\
Y_{t,x} (\left\{ s(\vr, \vc{m},E) \geq s_0 \right\}) = 1 \ \mbox{for a.a.} \ (t,x) .
\]

\subsection{Weak--strong uniqueness}

As shown in \cite[Theorem 3.3]{BreFei17},
the renormalized (DMV) solutions coincide with strong solutions emanating from the same initial data. The same result can be shown for the (DMV) solutions
in the sense of Definition \ref{D1}. The proof requires only obvious modification. In this context, the weak--strong uniqueness principle reads:

\Cbox{Cgrey}{

\begin{Theorem}{\bf[Weak (measure-valued) - strong uniqueness principle]} \label{T1}

Let the thermodynamic functions $p$, $e$, and $s$ satisfy Gibbs' equation (\ref{i4}), and the thermodynamic stability condition (\ref{i16}). In addition,
let the pressure be related to the internal energy through the caloric equation of state (\ref{i12bis}). Suppose that the Euler system (\ref{i7}--\ref{i10}) admits
a smooth ($C^1$) solution $(\vr, \vc{m}, E)$ in $[0,T)$ originating from the initial data $(\vr_0, \vm_0, E_0)$,
\[
\vr_0 \in C^1(\Omega), \ \vr_0 > 0, \ \vc{m}_0 \in C^1(R^N; R^N) , \ E_0 \in C^1(R),\ E_0 - \frac{1}{2} \frac{|\vc{m}_0|^2}{\vr_0} > \frac{\Ov{p}}{\gamma - 1} \vr_0^\gamma.
\]
Let $[\{ Y_{t,x} \}, \mathcal{D} ]$ be a (DMV) solution (renormalized (DMV) solution) in the sense of Definition \ref{D1} (Definition \ref{D2}) starting
from the same initial data,  meaning
\[
Y_{0,x} = \delta_{(\vr_0(x), \vc{m}_0(x), E_0(x)) } \ \mbox{for a.a.}\ x \in \Omega.
\]
Then $\mathcal{D} = 0$, and
\[
Y_{t,x} = \delta_{(\vr(t,x), \vc{m}(t,x), E(t,x)) } \ \mbox{for a.a.}\ (t,x) \in (0,T) \times \Omega.
\]

\end{Theorem}

}

\section{Measure--valued solutions generated in the vanishing dissipation limit}
\label{ZD}

The measure--valued solutions are natural candidates for describing the zero dissipation limits of more complex systems of Navier--Stokes type. Here, we show two results
in this direction. We restrict ourselves to the physically relevant case $N = 3$ although the same proof can be adapted to the case $N=1,2$.

\subsection{Vanishing dissipation limit of the Navier--Stokes--Fourier system}

The full compressible Navier--Stokes--Fourier system describing the motion of a general viscous and heat conducting fluid reads:

\begin{eqnarray}
\label{ZD1} \partial_t \vr + \Div (\vr \vu) &=& 0, \\
\label{ZD2} \partial_t (\vr \vu) + \Div \left(\vr {\vu \otimes \vu}\right) + \Grad p  &=& \Div \mathbb{S}, \\
\label{ZD3} \partial_t \left( \vr e \right) +
\Div \left( \vr e \vu \right) + \Grad \vc{q} &=& \mathbb{S} : \Grad \vu - p \Div \vu.
\end{eqnarray}

\noindent The viscous stress $\mathbb{S}$ and the heat flux $\vc{q}$ are determined by Stokes' law
\begin{equation} \label{ZD4}
\mathbb{S} = \mu \left( \Grad \vu + \Grad^t \vu - \frac{2}{3} \Div \vu \mathbb{I} \right) + \eta \Div \vu \mathbb{I},
\end{equation}
and Fourier's law
\begin{equation} \label{ZD5}
\vc{q} = - \kappa \Grad \vt.
\end{equation}
Finally, the relevant entropy balance is
\begin{equation} \label{ZD6}
\partial_t (\vr s) + \Div (\vr s \vu) + \Grad \left( \frac{\vc{q}}{\vt} \right) = \frac{1}{\vt} \left(
\mathbb{S} : \Grad \vu - \frac{\vc{q} \cdot \Grad \vt}{\vt} \right).
\end{equation}

\subsection{Weak solutions}

We adopt the concept of weak solution to the Navier--Stokes--Fourier system introduced in \cite[Chapters 2, 3]{FeNo6}.

\begin{Definition} \label{D3}

We say that $(\vr, \vt, \vu)$ is a \emph{weak solution} of the Navier--Stokes--Fourier system  if:
\begin{itemize}
\item $\vr \geq 0$, $\vt > 0$ a.a. in $(0,T) \times \Omega$;
\item
the equations (\ref{ZD1}), (\ref{ZD2}) are satisfied in the sense of (space periodic) distributions;
\item the entropy balance (\ref{ZD6}), is relaxed to the inequality
\[
\partial_t (\vr s(\vr, \vt)) + \Div (\vr s(\vr, \vt) \vu) + \Div \left( \frac{\vc{q}}{\vt} \right)
\geq \frac{1}{\vt} \left( \mathbb{S} : \Grad \vu - \frac{\vc{q} \cdot \Grad \vt}{\vt} \right)
\]
satisfied in the sense of distributions;
\item
the total energy balance
\bFormula{ZD8}
\intO{ \left[ \frac{1}{2} \vr |\vu|^2 + \vr e(\vr, \vt) \right] (\tau, \cdot) } =
\intO{ \left[ \frac{1}{2} \vr_0 |\vu_0|^2 + \vr_0 e(\vr_0, \vt_0) \right] }
\eF
holds for a.a. $\tau \in [0,T]$.
\end{itemize}

\end{Definition}

\begin{Remark} \label{Rws1}

Furthermore, a weak solution $(\vr, \vt, \vu)$ must belong to a certain regularity class for the weak formulation to make sense. The reader may consult
\cite[Chapters 2, 3]{FeNo6} for details.

\end{Remark}

\subsection{Constitutive equations}
\label{ce}

In analogy with Section \ref{CV}, we consider the pressure of the monoatomic gas related to the internal energy through
\begin{equation} \label{ZD9}
p = \frac{2}{3} \vr e, \ \mbox{meaning}\ \gamma = \frac{5}{3}.
\end{equation}
As shown in \cite[Chapter 2]{FeNo6}, relation (\ref{ZD9}) implies that
\begin{equation} \label{ZD10}
p(\vr, \vt) = \vt^{5/2} P \left( \frac{\vr}{\vt^{3/2}} \right)
\end{equation}
for some function $P$, cf. (\ref{i17}). In agreement with \cite[Chapter 3]{FeNo6}, we further assume that
 $P \in C^1 [0, \infty) \cap C^5(0, \infty)$ satisfies
\begin{equation} \label{ZD11}
P(0) = 0, \ P'(q) > 0 \ \mbox{for all} \ q \geq 0,
\end{equation}
\begin{equation} \label{ZD12}
0 < \frac{ \frac{5}{3} P(q) - P'(q) q }{q} < c \ \mbox{for all}\ q > 0, \ \lim_{q \to \infty} \frac{P(q)}{q^{5/3}} = \Ov{p} > 0.
\end{equation}
All the above requirements are just consequences of the thermodynamic stability hypothesis (\ref{i16}), except the stipulation
$\Ov{p} > 0$. Note that the standard pressure law $p = \vr \vt$ corresponds to $P(Z) = Z$, $\Ov{p} = 0$.

In agreement with (\ref{ZD9}), we set
\begin{equation} \label{ZD13}
e(\vr, \vt) = \frac{3}{2} \vt \left( \frac{ \vt^{3/2} }{\vr} \right) P \left( \frac{\vr}{\vt^{3/2}} \right),
\end{equation}
and, by virtue of Gibbs' relation (\ref{i4}),
\begin{equation} \label{ZD14}
s (\vr, \vt) = S \left( \frac{\vr}{\vt^{3/2}} \right),
\end{equation}
where
\begin{equation} \label{ZD15}
S'(q) = - \frac{3}{2} \frac{ \frac{5}{3} P(q) - P'(q) q }{q^2} < 0.
\end{equation}
We also impose the third law of thermodynamics in the form
\begin{equation} \label{ZD16}
\lim_{q \to \infty} S(q) = 0.
\end{equation}

\begin{Remark} \label{Rgib}

In the context of viscous fluids, it is convenient to work with the variables $(\vr, \vt, \vu)$. The function $S$ in (\ref{ZD14})
is therefore not the same as its counterpart expressed in the conservative variables in (\ref{entro}). Indeed the function $S = S(Z)$ 
in (\ref{entro}) is expressed in terms of $Z = p/\vr^\gamma$, while $S$ in (\ref{ZD14}) is a function of $q = \vr/ \vt^{c_v}$. If, for instance, $p = \vr \vt$, 
we obtain $q = Z^{1 - \gamma}$.

\end{Remark}

Finally, we suppose that the transport coefficients $\mu$, $\eta$ and $\kappa$ are continuously differentiable functions of the temperature $\vt$,
\begin{equation} \label{ZD17}
\mu, \eta \in C^1([0, \infty)),\
| \mu'(\vt) | \leq c, \ \underline{\mu} (1 + \vt) \leq \mu(\vt) ,\ 0 \leq \eta(\vt) \leq \Ov{\eta}(1 + \vt)
\  \mbox{for all}\ \vt \geq 0
\end{equation}
for certain constants $\underline{\mu} > 0$, $\Ov{\eta} > 0$ and
\begin{equation} \label{ZD18}
\kappa \in C^1([0, \infty)),\ \underline{\kappa} (1 + \vt^3) \leq \kappa(\vt) \leq \Ov{\kappa}(1 + \vt^3) \ \mbox{for all}\ \vt \geq 0
\end{equation}
for certain constants $\underline{\kappa} > 0$, $\Ov{\kappa} > 0$.

\subsection{Existence of weak solutions for the Navier--Stokes--Fourier system}

Motivated by the existence theory developed in \cite[Chapter 3]{FeNo6}, we consider the following system:
\begin{eqnarray}
\label{ZD19} \partial_t \vr + \Div (\vr \vu) &=& 0, \\
\label{ZD20} \partial_t (\vr \vu) + \Div \left(\vr {\vu \otimes \vu}\right) + \Grad (p + a p_R)  &=& \nu \Div \mathbb{S}, \\
\label{ZD21} \partial_t \left( \vr (e + a e_R) \right) +
\Div \left( \vr (e + a e_R) \vu \right) + \omega \Grad \vc{q} &=& \nu \mathbb{S} : \Grad \vu - p \Div \vu - \lambda (\vt - \Ov{\vt})^3.
\end{eqnarray}
Here $p_R$, $e_R$ are the radiation pressure and internal energy introduced in \cite[Chapter 3]{FeNo6},
\[
p_R = \frac{1}{3} \vt^4,\ e_R = \frac{\vt^4}{\vr}.
\]
The radiation components are multiplied by a (small) constant $a > 0$.
They provide a regularizing effect necessary for the existence theory developed in
\cite[Chapter 3]{FeNo6}.
The internal energy (\ref{ZD21}) contains a source term $- \lambda (\vt - \Ov{\vt})^3$, $\lambda > 0$ that may be interpreted as radiative ``cooling'' above a threshold
temperature $\Ov{\vt}$. The presence of this term provides a certain stabilizing effect necessary to perform the vanishing dissipation limit, cf. also \cite{FeMuNoPo}.
Note that the associated total energy and entropy balance read
\begin{equation} \label{ZD22b}
\left[ \intO{ \left( \frac{1}{2} \vr |\vu|^2 + \vr e + a \vt^4  \right) } \right]_{t = 0}^{t = \tau} +
\lambda \int_0^\tau \intO{ (\vt - \Ov{\vt})^3 } \dt = 0,
\end{equation}
\begin{equation} \label{ZD22}
\partial_t (\vr (s + a s_R) ) + \Div (\vr (s + a s_R) \vu) + \omega \Div \left( \frac{\vc{q}}{\vt} \right)
\geq \frac{1}{\vt} \left( \nu \mathbb{S} : \Grad \vu - \omega \frac{\vc{q} \cdot \Grad \vt}{\vt} \right) + \lambda \frac{(\Ov{\vt} - \vt)^3}{\vt}
\end{equation}
with
\[
s_R = \frac{4}{3} \frac{\vt^3}{\vr}.
\]
As stated in \cite[Chapter 3, Theorem 3.1]{FeNo6}, the problem (\ref{ZD19}--\ref{ZD21}) admits a global-in-time weak solution in the sense of Definition
\ref{D3}, whenever $a$, $\nu$, $\omega$, and $\lambda$ are positive and the constitutive restrictions specified in Section \ref{ce} hold.

\subsection{The asymptotic limit}

Our goal is to send $a \to 0$, $\nu \to 0$, $\omega \to 0$, and $\lambda \to 0$ to recover a dissipative measure--valued solution of the Euler system.
To this end, the following issues will be addressed:
\begin{itemize}
\item Uniform bounds based on the energy estimates that will guarantee boundedness of the state variables $(\vr, \vm, E)$.
\item Showing that the dissipation terms vanish in the asymptotic limit.
\item Identifying the dissipation defect $\mathcal{D}$ as well as the Young measure $\{ Y_{t,x} \}$ associated to the family of weak solutions.

\end{itemize}

\subsubsection{Uniform bounds}

The total energy balance (\ref{ZD22b}) yields immediately
\begin{equation} \label{ZD23}
{\rm ess} \sup_{t \in (0,T)} \intO{ \left[ \frac{|\vc{m}|^2}{\vr}  +  \vr e  + a \vt^4  \right]} +
\lambda \int_0^T \intO{ \vt^3 } \dt  \leq c(a, {\rm data}).
\end{equation}
In addition, it follows from hypotheses (\ref{ZD11}), (\ref{ZD12}), and (\ref{ZD13}) that
\[
\vr e(\vr, \vt) \ageq \vr \vt + \vr^{5/3};
\]
whence, in view of (\ref{ZD23}),
\begin{equation} \label{ZD24}
{\rm ess} \sup_{t \in (0,T)} \| \vr(t, \cdot) \|_{L^{5/3}(\Omega)} \leq c(a, {\rm data}),\
{\rm ess} \sup_{t \in (0,T)} \| \vr \vt (t, \cdot) \|_{L^{1}(\Omega)} \leq c(a, {\rm data}).
\end{equation}
Finally, writing $\vc{m} = \sqrt{\vr} \sqrt{\vr} \vu$, we deduce from (\ref{ZD23}), (\ref{ZD24}) that
\begin{equation} \label{ZD25}
{\rm ess} \sup_{t \in (0,T)} \| \vm (t, \cdot) \|_{L^{5/4}(\Omega)} \leq c(a, {\rm data}).
\end{equation}
Here and hereafter, the symbol $a \aleq b$ means $a \leq cb$ for a certain constant $c>0$.

Next, we have to handle the terms in the entropy balance. Writing
\[
\vr |s| \leq \frac{1}{2} \vr + \frac{1}{2} \vr s^2,\
\vr |s| \vu \leq \frac{1}{2} \vr |\vu|^2 + \frac{1}{2} \vr s^2
\]
we can see that it is enough to control $\vr s^2$ in $L^q$ for some $q > 1$. To this end, we use the third law of thermodynamics encoded in hypothesis
(\ref{ZD16}):
\[
|s| \aleq 1 \ \mbox{whenever} \ \vt^{3/2} \leq \vr.
\]
If $\vr < \vt^{3/2}$, we deduce from (\ref{ZD15}) that
\[
\vr s^2 \aleq \vr |\log(\vr)|^2 + \vr |\log(\vt)|^2,
\]
where $\vr \log(\vr)$ is controlled in the full range of $\vr$'s by (\ref{ZD24}). As for $\vr |\log(\vt)|^2$ it is dominated by $\vr \vt$ as long as $\vt \geq 1$. Thus it remains to control $\vr |\log(\vt)|^2$ in the range $\vr < \vt^{3/2}$, $\vt < 1$:
\[
\vr |\log(\vt)|^2 \leq \vr |\log( \vr^{2/3} )| = \frac{2}{3} \vr |\log(\vr)|.
\]
We may infer that
\begin{equation} \label{ZD26}
{\rm ess} \sup_{t \in (0,T)} \left\| \vr s \right\|_{L^q(\Omega)} +
{\rm ess} \sup_{t \in (0,T)} \left\| \vr s \vu \right\|_{L^q(\Omega; R^3)} \leq c(a, {\rm data})
\ \mbox{for some}\ q > 1.
\end{equation}

Note that all estimates obtained so far are uniform with respect to the parameters $\nu$, $\omega$, $\lambda$, and $a$ as long as $a \aleq 1$.
They are strong enough to pass to the limit in the system (\ref{ZD19}--\ref{ZD21}) to generate a (DMV) solution of the limit Euler system as soon as
we show that the dissipative terms vanish in the asymptotic regime. First observe that (\ref{ZD22}), together with
hypotheses (\ref{ZD17}), (\ref{ZD18}), gives rise to the bound
\[
\nu \int_0^T \intO{ \left| \Grad \vu + \Grad \vu^t - \frac{2}{3} \Div \vu \mathbb{I} \right|^2 } \dt \leq c(a, {\rm data}),
\]
which, after a simple by parts integration, yields
\begin{equation} \label{ZD27}
\nu \int_0^T \intO{ |\Grad \vu |^2 } \leq c(a, {\rm data}).
\end{equation}
Seeing that the total mass of the fluid is conserved,
\begin{equation} \label{ZD28}
\intO{ \vr(t, \cdot) } = \intO{ \vr_0 } \ \mbox{for any}\ t \geq 0,
\end{equation}
we may use a version of Poincare's inequality to deduce from (\ref{ZD24}), (\ref{ZD27}), (\ref{ZD28}) that
\[
\nu \int_0^T \intO{ |\vu|^2 } \dt \aleq \nu \left[ \int_0^T \intO{ \vr |\vu|^2 }\dt + \int_0^T \intO{ |\Grad \vu |^2 } \dt \right] \leq c(a, {\rm data});
\]
whence
\begin{equation} \label{ZD29}
\nu \int_0^T \| \vu \|^2_{W^{1,2}(\Omega; R^3)} \dt \leq c(a, {\rm data}).
\end{equation}

Applying the same treatment to $\vt$, we get
\begin{equation} \label{ZD30}
\omega \int_0^T \intO{ \left( \frac{1}{\vt^2} + \vt \right) |\Grad \vt|^2 } \dt \aleq c(a, {\rm data}),
\end{equation}
and, using (\ref{ZD24}),
\[
\omega \int_0^T \intO{ |\vt|^2 } \dt \aleq \omega \left[ \int_0^T \left( \intO{ \vr \vt } \right)^2 \dt + \int_0^T \intO{ |\Grad \vt |^2 } \dt \right] \leq c(a, {\rm data});
\]
whence
\begin{equation} \label{ZD31}
\omega \int_0^T \| \vt \|^2_{W^{1,2}(\Omega)} \dt \leq c(a, {\rm data}).
\end{equation}

\subsubsection{Vanishing dissipation limit}

For the sake of simplicity, we suppose that
\[
\nu = \omega = \ep,\ a = \ep^\alpha,\ \lambda = \ep^\beta
\]
for suitable $\alpha > 0$, $\beta > 0$ fixed below. For $(\vre, \vte, \vue)_{\ep > 0}$ - a sequence of weak solutions of problem (\ref{ZD19}--\ref{ZD21}) in the sense of Definition \ref{D3} - we set
\[
(\vre, \vme, \Ee)_{\ep > 0} ,\ \vme = \vre \vue,\ \Ee = \frac{1}{2} \vre |\vue|^2 + \vre e(\vre, \vte).
\]
In addition,
we suppose that the initial data
\[
\left( \vr_{0,\ep}, \vr_{0,\ep} \vu_{0,\ep}, \frac{1}{2} \vr_{0,\ep} |\vu_{0,\ep} |^2 + \vr_{0,\ep} e(\vr_{0,\ep}, \vt_{0,\ep}) \right)_{\ep > 0}
\]
generate a Young measure $Y_{0,x}$, specifically,
\begin{equation} \label{ZD33}
\begin{split}
\intO{ \vr_{0,\ep} \phi } &\to \intO{ \left< Y_{0,x}; \vr \right> \phi } \ \mbox{for any}\ \phi \in \DC(\Omega); \\
\intO{ \vr_{0,\ep} \vu_{0,\ep} \cdot \boldsymbol{\phi} } &\to \intO{ \left< Y_{0,x}; \vm \right> \cdot \boldsymbol{\phi} } \ \mbox{for any}\
\boldsymbol{\phi} \in \DC(\Omega; R^3),\\
\intO{ \left[ \frac{1}{2} \vr_{0,\ep} |\vu_{0,\ep}|^2 + \vr_{0,\ep} e(\vr_{0,\ep},\vt_{0,\ep}) + \ep^\alpha \vt_{0,\ep}^4 \right] \phi }
&\to \intO{\left< Y_{0,x}; E \right> \phi } \ \mbox{for any} \ \phi \in \DC(\Omega);\\
\intO{ \left[ \vr_{0,\ep} s(\vr_{0,\ep}, \vt_{0,\ep}) + \ep^\alpha \frac{4}{3} \vt^3_{0,\ep} \right] \phi }
&\to \intO{ \left< Y_{0,x}; \vr s(\vr, \vc{m}, E) \right> \phi } \ \mbox{for any}\ \phi \in \DC(\Omega).
\end{split}
\end{equation}

In view of the uniform bounds (\ref{ZD23}--\ref{ZD25}) and the fundamental
theorem of the theory of Young measures (see e.g. Ball \cite{BALL2}), there is a subsequence of $(\vre, \vme, \Ee)_{\ep > 0}$ (not relabeled here) that generates
a Young measure $\{ Y_{t,x} \}_{(t,x) \in (0,T) \times \Omega}$. Moreover, passing to the limit in the total energy balance (\ref{ZD22b}), we obtain
\begin{equation} \label{ZD34}
\left[ \intO{ \left< Y_{\tau,x} ; E \right> } \right]_{t = 0}^\tau + \mathcal{D} (\tau) = 0;
\end{equation}
for a.a. $\tau \in (0,T)$, where
\begin{equation} \label{ZD35}
\mathcal{D}(\tau) \geq \liminf_{\ep \to 0} \intO{ \left[ \frac{1}{2} \vre |\vue|^2 +
\vre e(\vre, \vte) + \ep^\alpha \vte^4 \right] } - \intO{ \left< Y_{\tau,x} ; E \right> } \ \mbox{for a.a}\ \tau \in (0,T).
\end{equation}
Indeed the extra term in (\ref{ZD22b}) can be handled as
\begin{equation} \label{ZD36}
\lambda \int_0^\tau \intO{ (\vte - \Ov{\vt})^3 } \ageq - \ep^\beta \int_0^\tau \intO{ \left[ \vte^2 \Ov{\vt} + \Ov{\vt}^2 \vte + \Ov{\vt^3} \right] } \dt
\approx - \ep^{\beta/3} \to 0,
\end{equation}
where we have used the bound (\ref{ZD23}).
Moreover, it is easy to pass to the limit in the weak formulation of (\ref{ZD19}), to obtain (\ref{MV1}).

Next, seeing that
\[
\nu \mathbb{S}(\vte, \Grad \vue) \approx \ep \vte |\Grad \vue| = \ep^{1/2} \vte \ep^{1/2} |\Grad \vue|,
\]
we may use the bounds (\ref{ZD23}), (\ref{ZD27}) to conclude that
\[
\left\| \nu \mathbb{S}(\vte, \Grad \vue) \right\|_{L^1((0,T) \times \Omega)} \aleq \ep^{\frac{1}{2} - \frac{\beta}{3}}.
\]
Thus choosing $0 < \beta < \frac{3}{2}$ we can pass to the limit in the momentum balance (\ref{ZD20}) to recover
(\ref{MV2}). Note that the measure $\mu_C$ contains the concentration defect of the terms
\[
\vre \vue \otimes \vue,\ p(\vre, \vue), \ \mbox{and}\ \ep^\alpha \vte^4
\]
and, by virtue of (\ref{ZD35}), it is controlled by $\mathcal{D}$ exactly as required in (\ref{MV5}).

Finally, it remains to perform the limit in the entropy balance
(\ref{ZD22}) to obtain (\ref{MV4}). First, by virtue of the same argument as in (\ref{ZD36}), we get
\[
\lambda \frac{ (\Ov{\vt} - \vte)^3 }{\vte} \ageq - \ep^{\beta/3} \to 0.
\]

Next, the entropy heat flux can be treated as
\[
\omega \frac{\vc{q(\vte)}}{\vte} = \ep \frac{\kappa_1(\vte)}{\vte} \Grad \vte + \ep \frac{\kappa_2(\vte)}{\vte} \Grad \vte,
\]
where $\kappa_1$ is bounded and $\kappa_2 \approx \vt^3$. In view of (\ref{ZD30}), we get
\[
\ep \left| \frac{\kappa_1(\vte)}{\vte} \Grad \vte \right|_{L^2((0,T) \times \Omega)} \aleq \ep \left\| \  |\Grad \log(\vte)| + |\Grad \vte | \
\right\|_{L^2((0,T) \times \Omega)} \aleq \ep^{1/2}.
\]
Moreover,
\[
\ep \vte^2 \Grad \vte = \ep \frac{1}{3} \Grad \vte^3,
\]
where, in accordance with (\ref{ZD23})
\[
\ep \vte^3 \to 0 \ \mbox{in}\ L^1((0,T) \times \Omega) \ \mbox{as soon as}\ 0 < \beta < 1.
\]

Seeing that the terms $\vre s(\vre, \vte)$ and $\vre s(\vre, \vte) \vue$ are controlled via (\ref{ZD26}), it remains to show
\[
a \left( {\rm ess} \sup_{t \in (0,T)} \| \vte^3 \|_{L^1(\Omega)} + \left\| \vte^3 \vue \right\|_{L^1((0,T) \times \Omega)} \right) \to 0.
\]
Seeing that, in view of the energy bound (\ref{ZD23}),
\[
 {\rm ess} \sup_{t \in (0,T)} a\| \vte^3 \|_{L^1(\Omega)} \to 0,
\]
we have to handle only the second term.
To this end, write
\[
a \vte^3 \vue = \ep^\alpha \nu^{-1/2} \vte^3 \nu^{1/2} \vue.
\]
In view of (\ref{ZD29}), it is enough to show that
\[
\ep^{\alpha - \frac{1}{2} } {\rm ess}\sup_{t \in (0,T)} \| \vte^3 \|_{L^{4/3}(\Omega)} \to 0.
\]
However, this follows from the energy bound (\ref{ZD23}) as soon as $\alpha > 2$.

We have shown the following result:

\Cbox{Cgrey}{

\begin{Theorem} \label{T2}

Suppose that the thermodynamic functions $p$, $e$, and $s$ satisfy (\ref{ZD9}--\ref{ZD16}), and the transport coefficients
$\mu$, $\eta$ and $\kappa$ satisfy (\ref{ZD17}), (\ref{ZD18}). Let $(\vre, \vte, \vue)_{\ep > 0}$ be a family of weak solutions
to the Navier--Stokes--Fourier system (\ref{ZD19}--\ref{ZD21}), where
\[
\nu = \omega = \ep, \ a= \ep^\alpha,\ \alpha > 2,\ \lambda = \ep^\beta,\ 0 < \beta < 1.
\]
Finally, let the initial data $(\vr_{0, \ep}, \vt_{0, \ep}, \vu_{0, \ep} )_{\ep > 0}$ generate a Young measure $Y_{0,x}$ specified in (\ref{ZD33}).

Then (at least for a suitable subsequence)
\[
\left( \vre, \vre \vue, \frac{1}{2} \vre |\vue|^2 + \vre e(\vre, \vte) \right)_{\ep > 0}
\]
generates a Young measure $\{ Y_{t,x} \}_{(t,x) \in (0,T) \times \Omega}$ and a dissipation defect $\mathcal{D}$ specified in
(\ref{ZD35}) that represent a dissipative measure--valued solution of the Euler system (\ref{i7}--\ref{i10})
in the sense of Definition \ref{D1}.

\end{Theorem}

}

In view of Theorem \ref{T1}, we immediately obtain the following corollary that can be seen as a version of the result in \cite{Fei2015A}.

\begin{Corollary} \label{C1}

In addition to the hypotheses of Theorem \ref{T2} suppose that the limit Euler system (\ref{i7}--\ref{i9}) admits a smooth ($C^1$) solution $(\vr, \vm, E)_{\ep > 0}$ in
$[0,T] \times \Omega$.

Then
\[
\vre \to \vr,\ \vre \vue \to \vc{m}, \frac{1}{2} \vre |\vue|^2 + \vre e(\vre, \vte) \to E \ \mbox{in}\ L^1((0,T) \times \Omega).
\]

\end{Corollary}

Indeed the fact that the limit DMV solution is represented by the Dirac masses implies (up to a subsequence) strong a.a. pointwise convergence. In addition,
the limit defect $\mathcal{D}$ vanishes which implies strong convergence in the $L^1-$norm.

\section{Vanishing dissipation limit of Brenner's model}
\label{BM}

The measur--valued solutions constructed via the Navier--Stokes--Fourier system do not, or at least are not known, to satisfy the renormalized
entropy balance (\ref{MV6}). In addition, the hypotheses imposed in Theorem \ref{T2} on the constitutive relations and the transport coefficients are rather awkward and, unfortunately, do not cover the mostly used Boyle--Mariotte law $p = \vr \vt$, with $e = c_v \vt$.
It is interesting to see that a model proposed by H.Brenner \cite{BREN2}, \cite{BREN}, \cite{BREN1}  behaves actually better in the vanishing dissipation limit and does not suffer the above mentioned
drawbacks.

Brenner's model consists of the following equations:

\Cbox{Cgrey}{

\begin{eqnarray}
\label{b1} \partial_t \vr + \Div (\vr \vv_m)&=& 0, \\
\label{b2} \partial_t (\vr \vv) + \Div \left(\vr \vv\otimes \vv_m\right) + \Grad p(\vr,\vt) &=& {\ep} \Div \mathbb{S} ,
\\
\begin{split}
\label{b3} \partial_t \left(\vr(\frac12|\vv|^2 +  e(\vr, \vt))\right) +
\Div \left(\vr(\frac12|\vv|^2 + e(\vr, \vt))\vv_m \right)\\ + \Div(p(\vr,\vt)\vv) + {\ep}\Div \vq &=& {\ep}\Div(\mathbb S \vv),
\end{split}
\end{eqnarray}
}

\noindent
where, similarly to the preceding section, $\ep >0$ is a small parameter supposed to vanish in the zero dissipation limit.
Brenner's main idea was to introduce two velocity fields - $\vv$ and $\vv_m$ - interrelated through
\begin{equation}\label{b4}
\vv - \vv_m = {\ep} K \Grad \log (\vr),
\end{equation}
where $K \geq 0$ is a purely {\it phenomenological} coefficient. Similarly to the previous part, we suppose
\begin{equation}\label{b8}
\vq = -\kappa \Grad \vt,
\end{equation}
and
\begin{equation}
\label{b10}
\S = \mu\left(\Grad \vv + \Grad^T \vv - \frac23 \Div \vv \I\right)
 + \eta \Div \vv \I.
\end{equation}

As a matter of fact, Brenner's model has been thoroughly criticized and its relevance to fluid mechanics questioned in Oettinger et al. \cite{OeStLi}. On the other hand,
it is mathematically tractable and yields essentially better theory than the standard Navier--Stokes--Fourier system, see e.g. \cite{FeiVas}, Cai, Cao, Sun \cite{CaiCaoSun}.
Recently, the interest in ``two velocity models'' has been revived in Bresch et al. \cite{BrDeZa}, \cite{BrGiZa}.

Leaving apart the conceptual difficulties of the model, we claim that it generates in the vanishing dissipation limit a \emph{renormalized} (DMV) solution of the Euler
system (\ref{i7}--\ref{i10}). The crucial aspect of the analysis is a specific form of the coefficient $K$
in (\ref{b4}). Note that $K$ is taken constant in \cite{FeiVas} as well as in  Cai et al. \cite{CaiCaoSun}, while
Brenner proposed $K = \frac{\kappa}{c_p \vr}$, see \cite{BREN}, where $c_p$ denotes the specific heat at constant \emph{pressure}.
Further, we assume Boyle--Mariotte law
\begin{equation} \label{b102}
p = \vr \vt,\ e = c_v \vt, \ s = \log(\vt^{c_v}) - \log(\vr),
\end{equation}
where $c_v = \frac{1}{\gamma - 1}$ is the specific heat at constant \emph{volume}.
As we show below, a convenient form of $K$ reads
\begin{equation} \label{b100}
K = \frac{\kappa}{c_v \vr}.
\end{equation}

\subsection{Estimates of entropy}

To begin, let us say frankly that there is no rigorous existence result for the Brenner model if $K$ is given by (\ref{b100}). For this reason, we restrict ourselves
to showing the argument how the renormalized entropy inequality is obtained. Furthermore, we shall suppose that we deal with a smooth solution of (\ref{b1}--\ref{b3}).

Taking the scalar product of the momentum equation (\ref{b2}) on $\vv$ and subtracting the resulting expression from (\ref{b3}), we obtain the internal energy balance
\begin{equation} \label{b101}
\partial_t (\vr e(\vr, \vt)) + \Div (\vr e (\vr, \vt) \vv_m )- \ep \Div \left(\kappa \Grad \vt \right) = \ep \mathbb{S} : \Grad \vu - p(\vr, \vt) \Div \vv.
\end{equation}
Dividing (\ref{b101}) on $\vt$ and keeping in mind the constitutive relations (\ref{b102}), (\ref{b100}), we obtain
\[
\partial_t (\vr \log (\vt^{c_v})) + \Div \left( \vr \log( \vt^{c_v}) \vv_m \right) - \ep \Div \left( \frac{\kappa} {\vt} \Grad \vt \right) =
\frac{1}{\vt} \left( \mathbb{S} : \Grad \vv + \frac{\kappa}{\vt} |\Grad \vt|^2 \right) - \vr \Div \vv.
\]
Furthermore, in accordance with (\ref{b4}),
\[
\vr \Div \vv = \vr \Div (\vv - \vv_m) + \vr \Div \vv_m =
\ep \vr \Div \left( \frac{\kappa}{\vr c_v} \Grad \log(\vr) \right) - \partial_t (\vr \log(\vr)) - \Div \left( \vr \log(\vr) \vv_m \right).
\]
Thus we have deduced the entropy balance
\[
\begin{split}
\partial_t (\vr s(\vr, \vt) ) &+ \Div \left( \vr s(\vr, \vt) \vv_m \right) -
\ep \Div \left[ \frac{\kappa}{c_v} \Grad \left( \log(\vt^{c_v}) -  \log(\vr) \right) \right]  \\
&= \frac{1}{\vt} \left( \mathbb{S} : \Grad \vv + \frac{\kappa}{\vt} |\Grad \vt|^2 \right) + \ep \frac{\kappa}{c_v} |\Grad \log(\vr) |^2,
\end{split}
\]
or, in a slightly different form,
\begin{equation} \label{b103}
\begin{split}
\partial_t (\vr s(\vr, \vt) ) &+ \Div \left( \vr s(\vr, \vt) \vv_m \right) -
\ep \Div \left[ \frac{\kappa}{c_v} \Grad s(\vr, \vt) \right]  \\
&= \frac{1}{\vt} \mathbb{S} : \Grad \vv +  {\kappa} |\Grad \log{\vt}|^2  + \ep \frac{\kappa}{c_v} |\Grad \log(\vr) |^2.
\end{split}
\end{equation}

Multiplying (\ref{b103}) on $Z'(s)$, where $Z$ is a non-decreasing concave function, we obtain the desired relation
\begin{equation} \label{b104}
\begin{split}
\partial_t &(\vr Z(s(\vr, \vt)) ) + \Div \left( \vr Z(s(\vr, \vt)) \vv_m \right) -
\ep \Div \left[ \frac{\kappa}{c_v} \Grad Z(s(\vr, \vt)) \right]  \\
&= \frac{Z'(s)}{\vt} \mathbb{S} : \Grad \vv +  {\kappa} Z'(s) |\Grad \log{\vt}|^2  + \ep Z'(s) \frac{\kappa}{c_v} |\Grad \log(\vr) |^2
- Z''(s) \frac{\kappa}{c_v} |\Grad s(\vr, \vt) |^2.
\end{split}
\end{equation}
In particular, by virtue of the parabolic maximum principle, we may deduce from (\ref{b104}) that
\[
\vr^{\gamma} \leq c(s_0) c_v \vr \vt,
\]
where $s_0 > - \infty$ is the initial entropy, see \cite[Section 2.1.1]{BreFei17}.

Similarly to \cite[Section 2.1]{BreFei17} we could show that a family of (strong) solutions $(\vre, \vte, \vue)_{\ep > 0}$
of Brenner's model generates a measure--valued solution
of the Euler system (\ref{i7}--\ref{i10}) satisfying (\ref{MV6}) in the asymptotic limit $\ep \to 0$. As the result is only formal (we do not know whether the strong or even weak solutions exist) we do not state it as a theorem and leave the interested reader to work out the details.

\def\cprime{$'$} \def\ocirc#1{\ifmmode\setbox0=\hbox{$#1$}\dimen0=\ht0
  \advance\dimen0 by1pt\rlap{\hbox to\wd0{\hss\raise\dimen0
  \hbox{\hskip.2em$\scriptscriptstyle\circ$}\hss}}#1\else {\accent"17 #1}\fi}


\end{document}